\documentclass[12pt]{article}

\usepackage{pstricks,pst-node}

\usepackage{amsmath,amsthm}
\usepackage{amssymb}

\usepackage{color}
\usepackage{xspace}
\usepackage[colorlinks=true,
linkcolor=brown,
filecolor=brown,
citecolor=red]{hyperref}

\def\mbf#1{\mathchoice{\hbox{\boldmath $\displaystyle #1$}}
        {\hbox{\boldmath $\textstyle #1$}}
        {\hbox{\boldmath $\scriptstyle #1$}}
        {\hbox{\boldmath $\scriptscriptstyle #1$}}} % mathboldface

\begin{document}
\newtheorem{theorem}{Theorem}
\newtheorem{defn}[theorem]{Definition}
\newtheorem{lemma}[theorem]{Lemma}
\newtheorem{prop}[theorem]{Proposition}
\newtheorem{cor}[theorem]{Corollary}
%\vspace*{-5mm}
\begin{center}
{\Large
A sign-reversing involution to count labeled \\ 
lone-child-avoiding trees    \\ 
}

\vspace{5mm}
David Callan  \\
\noindent {\small Department of Statistics, 
University of Wisconsin-Madison,  Madison, WI \ 53706}  \\
{\bf callan@stat.wisc.edu}  \\
June 24 2014
\end{center}

\begin{abstract}
We use a sign-reversing involution to show that trees on the vertex set $[n]$, considered to be rooted at 1, in which no vertex has exactly one child are counted by $\frac{1}{n}\sum_{k=1}^{n} (-1)^{n-k}\binom{n}{k}\frac{(n-1)!}{(k-1)!} k^{k-1}$. This result corrects a persistent misprint in the Encyclopedia of Integer Sequences. %A108919

\end{abstract}

\vspace{2mm}

\section{Introduction} \label{intro}

\vspace*{-4mm}

A graph is said to be \emph{series-reduced} or \emph{homeomorphically irreducible} if no vertex has degree 2, and a tree is a connected graph with no cycles. 
Several kinds of series-reduced trees have been enumerated, such as 
\emph{free} (sequence \htmladdnormallink{A000014}{http://oeis.org/A000014} in OEIS \cite{oeis}\,),
\emph{rooted} (\htmladdnormallink{A059123}{http://oeis.org/A059123}),
\emph{labeled} (\htmladdnormallink{A005512}{http://oeis.org/A005512}), and 
\emph{rooted labeled} (\htmladdnormallink{A060313}{http://oeis.org/A060313}),
the latter two differing only by a factor of the number of vertices.
A \emph{planted} tree is a rooted tree in which the root has degree 1. 
Series-reduced planted trees are counted by 
\htmladdnormallink{A001678}{http://oeis.org/A001678}. 
Let us say a rooted tree is \emph{lone-child-avoiding} if no vertex has exactly one child. 
A series-reduced planted labeled tree on vertex set $\{0,1,2,\dots,n\}$ with root 0 is equivalent to 
a lone-child-avoiding tree on $[n]:=\{1,2,\dots,n\}
$: delete the root and its one edge and re-root the new tree at the child of the old root.

The number of lone-child-avoiding 
rooted trees on $[n] $ is, of course, $n$ times the number of such trees with root 1. 
We will show that the latter are counted by the alternating sum
\begin{equation}\label{eq1}
\frac{1}{n}\sum_{k=1}^{n} (-1)^{n-k}\binom{n}{k}\frac{(n-1)!}{(k-1)!} k^{k-1}\, ,
\end{equation} 
correcting an erroneous entry in OEIS that has stood for 40+ years, namely
\htmladdnormallink{A002792}{http://oeis.org/A002792}, which is
 1, 0, 1, 1, 13, 51, 601, 4806, 39173, 775351 (that's all)  and surely should be  the sequence
1, 0, 1, 1, 13, 51, 601, 4803, 63673, 775351, $\dots$ 
(\htmladdnormallink{A108919}{http://oeis.org/A108919}) generated by (\ref{eq1}).

While the standard methods of recurrence relations and generating functions 
can readily be applied to our problem (see \cite{meir}, \cite[Ex. 3.3.26]{gj} for similar applications), we will use a sign-reversing involution.  
Section \ref{configs} introduces weighted objects whose total weight is obviously given by (\ref{eq1}). 
Section \ref{invol} then presents a sign-reversing involution to show that the total weight is also equal to the number of trees being counted.

\section[G-configurations]{\protect{$\mbf{G}$}-configurations} \label{configs} 

\vspace*{-4mm}

A rooted $G$-configuration of size $n$ is a hybrid graph object on the vertex set $[n]$ with $n-1$ edges, some of which are directed and some not as follows. The undirected edges 
form a tree on some nonempty subset $V$ of $[n]$ with some element $r$ of $V$  
designated as the root of the tree. There is exactly one directed edge (arc) from 
each vertex in $[n]\, \backslash \,V$ to $[n]\, \backslash \,\{r\}$ (thus no arc goes into the root)
and no two arcs end at the same vertex. 
We distinguish between \emph{tree} vertices (in the tree) and \emph{arc} vertices (that start an arc). 
Figure 1 shows an example.

\begin{center} 
\begin{pspicture}(-8,-2.1)(5,3.7)%\showgrid
%\psset{unit=.8cm}
\psdots(-2,2)(0,3)(2,2)(1,1)(2,1)(3,1)(2,0)
{\psset{linewidth=1.5pt} 
\psdots(0,3)}
\psline(-2,2)(0,3)(2,2)(1,1)
\psline(2,0)(2,1)(2,2)(3,1)

\psbezier[linewidth=1pt,linecolor=blue,showpoints=false,arrowsize=1.5pt 2.5]{<-}(-7.1,3.1)(-7.4,3.6)(-6.6,3.6)(-6.9,3.12)

\psbezier[linewidth=1pt,linecolor=blue,showpoints=false,arrowsize=1.5pt 2.5]{<-}(-5.4,3.1)(-5.0,3.4)(-4.5,3.4)(-4.1,3.1)

\psbezier[linewidth=1pt,linecolor=blue,showpoints=false,arrowsize=1.5pt 2.5]{->}(-5.4,2.9)(-5.0,2.6)(-4.5,2.6)(-4.1,2.9)

\rput(-7,2.7){\textrm{{\blue \small $8$}}}
\rput(-5.7,3){\textrm{{\blue \small $9$}}}
\rput(-3.7,3){\textrm{{\blue \small $14$}}}

\rput(-3.8,2){\textrm{{\blue \small $5$}}}
\rput(-2,-0.3){\textrm{{\blue \small $2$}}}
\rput(-0.5,-0.3){\textrm{{\blue \small $6$}}}
\rput(1,-0.3){\textrm{{\blue \small $12$}}}

\psline[linecolor=blue,arrows=->,arrowsize=3pt 3.5](-3.4,2)(-2.2,2)

\psline[linecolor=blue,arrows=->,arrowsize=3pt 3.5](-1.9,0)(-0.6,0)
\psline[linecolor=blue,arrows=->,arrowsize=3pt 3.5](-0.4,0)(0.9,0)
\psline[linecolor=blue,arrows=->,arrowsize=3pt 3.5](1.1,0.1)(1.9,.9)

\rput(0,3.3){\textrm{{\small $1$}}}
\rput(-1.9,1.7){\textrm{{\small $3$}}}
\rput(2.2,2.3){\textrm{{\small $10$}}}
\rput(1,0.7){\textrm{{\small $4$}}}
\rput(3,0.7){\textrm{{\small $11$}}}
\rput(2.2,0.7){\textrm{{\small $7$}}}
\rput(2,-0.3){\textrm{{\small $13$}}}

{\psset{linecolor=blue} 
\psdots(-2,0)(-0.5,0)(1,0)
\psdots(-7,3)(-5.5,3)(-4,3)(-3.5,2)}

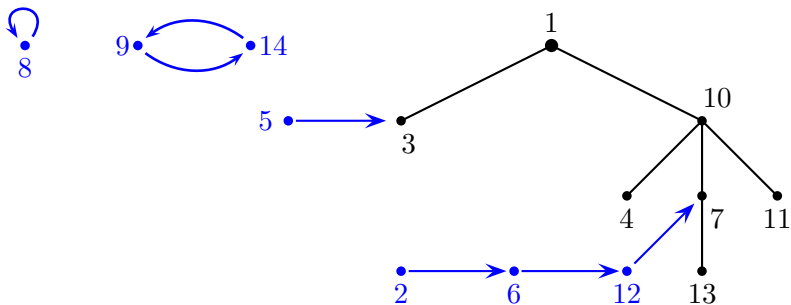
\rput(-1.5,-1.5){\textrm{{\small A G-configuration ( = rooted $G$-configuration with root 1) of size 14}}}
\rput(-1.5,-2.1){\textrm{{\small {\blue arcs and arc vertices in blue}}}}

\end{pspicture}
\end{center} 
\centerline{Figure 1}

\vspace*{2mm}

The number of G-configurations of size $n$ with $k$ tree vertices is $\binom{n}{k}$ [\,choose tree vertices\,] 
$\times\, k^{k-1}$ [\,form rooted tree, Cayley's formula \cite[Chap.\,30]{bookproofs}\,] $\times\, \underbrace{(n-1)(n-2)\  \cdots\ }_{n-k}$ 
[\,choose arcs\,]. 

Now we restrict the root to be 1. A (plain) $G$-configuration is a rooted 
$G$-configuration whose root is 1, dividing the count by a factor of $n$. Assign a weight of $(-1)^{n-k}$ 
to each $G$-configuration of size $n$ with $k$ tree vertices. Thus the total weight of all 
$G$-configurations of size $n$ is given by (\ref{eq1}).

\section{A sign-reversing involution} \label{invol}

\vspace*{-4mm}

We now define a sign-reversing involution on all $G$-configurations of size $n$ except the ones 
whose tree has $n$ vertices and no lone children. The involution converts a lone child to an 
arc vertex or vice versa, thereby changing the sign of the weight. The vertex $m$ to be 
converted is the maximum among all the arc vertices and the lone-child vertices. 
This maximum exists unless all vertices are in the tree and none is a lone child.

Suppose first that $m$ is an arc vertex. Either $m$ is a vertex in a cycle comprised of arcs or there is a path from $m$ along arcs terminating at some tree vertex $c$. In both cases, remove the arc that starts at $m$.
Then, in the first case, insert $m$ into the tree as the lone child of the root, 
redirecting into $m$ all the edges that originally 
went into the root (Figure 2a). An original arc $m\to m$ is lost, but an arc $a\to m$ with $a\ne m$ is retained. 

\begin{center} 
\begin{pspicture}(-8,-0.4)(8,3.3)%\showgrid

{\psset{linewidth=1.5pt} 
\psdots(-2.5,3)
\psdots(4.5,3)}

\psdots(-3,1)(-3,2)(-2,1)(-2,2)

\psline(-3,1)(-3,2)(-2.5,3)(-2,2)(-2,1)

\psbezier[linewidth=1pt,linecolor=blue,showpoints=false,arrowsize=1.5pt 2.5]{->}(-4,3.1)(-4.5,3.6)(-5.5,3.6)(-6,3.1)

\rput(-3,.7){\textrm{{\small $4$}}}
\rput(-2.75,1.9){\textrm{{\small $7$}}}
\rput(-1.75,1.9){\textrm{{\small $9$}}}
\rput(-2,.7){\textrm{{\small $5$}}}
\rput(-2.5,3.3){\textrm{{\small $1$}}}
 
\psline[linecolor=blue,arrows=->,arrowsize=3pt 3.5](-5.9,3)(-5.1,3)
\psline[linecolor=blue,arrows=->,arrowsize=3pt 3.5](-4.9,3)(-4.1,3)
\psline[linecolor=blue,arrows=->,arrowsize=3pt 3.5](-3.9,2)(-3.1,2)

\rput(-6,2.7){\textrm{{\blue \small $3$}}}
\rput(-5,2.7){\textrm{{\blue \small $2$}}}
\rput(-4,2.7){\textrm{{\blue \small $8$}}}
\rput(-4,1.7){\textrm{{\blue \small $6$}}}

\psdots(4,0)(4,1)(4.5,2)(5,1)(5,0)
\psline(4,0)(4,1)(4.5,2)(5,1)(5,0)
 \psline(4.5,2)(4.5,3)

{\psset{linecolor=blue} 
\psdots(2,2)(3,2)(3,1)
\psdots(-6,3)(-5,3)(-4,3)(-4,2) }

\psline[linecolor=blue,arrows=->,arrowsize=3pt 3.5](2.1,2)(2.9,2)
\psline[linecolor=blue,arrows=->,arrowsize=3pt 3.5](3.1,2)(4.4,2)

\psline[linecolor=blue,arrows=->,arrowsize=3pt 3.5](3.1,1)(3.9,1)

\rput(2,1.7){\textrm{{\blue \small $3$}}}
\rput(3,1.7){\textrm{{\blue \small $2$}}}
\rput(3,0.7){\textrm{{\blue \small $6$}}}

\rput(0.0,1){\textrm{$\rightarrow$}}

\rput(4,-.3){\textrm{{\small $4$}}}
\rput(4.25,.9){\textrm{{\small $7$}}}
\rput(5.25,.9){\textrm{{\small $9$}}}
\rput(5,-.3){\textrm{{\small $5$}}}
\rput(4.5,3.3){\textrm{{\small $1$}}}
\rput(4.8,2.1){\textrm{{\small $8$}}}

\rput(-4.5,-0.2){\textrm{{\small $9$ is not a lone child,}}}
\rput(-4.5,-0.7){\textrm{{\small and so $m=8$}}}

\end{pspicture}
\end{center} 
\centerline{Figure 2a}

\vspace*{2mm}

In the second case, insert $m$ into the tree as the lone child of $c$, 
redirecting into $m$ all the original child edges of $c$ (Figure 2b).

\begin{center} 
\begin{pspicture}(-8,-1.9)(8,3.3)%\showgrid
%\psset{unit=.8cm}

{\psset{linewidth=1.5pt} 
\psdots(-1.5,3)
\psdots(4.5,3)}

\psdots(-1.5,3)(-1.5,2)(-2,1)(-1,1)(-1,0)

\psline(-1.5,3)(-1.5,2)(-2,1)
\psline(-1.5,2)(-1,1)(-1,0)

\rput(-6,2.3){\textrm{{\blue \small $6$}}}
\rput(-5,2.3){\textrm{{\blue \small $9$}}}
\rput(-4,2.3){\textrm{{\blue \small $3$}}}
\rput(-3,2.3){\textrm{{\blue \small $4$}}}
\rput(-3,1.3){\textrm{{\blue \small $8$}}}
 
\psline[linecolor=blue,arrows=->,arrowsize=3pt 3.5](-5.9,2)(-5.1,2)
\psline[linecolor=blue,arrows=->,arrowsize=3pt 3.5](-4.9,2)(-4.1,2)
\psline[linecolor=blue,arrows=->,arrowsize=3pt 3.5](-3.9,2)(-3.1,2)
\psline[linecolor=blue,arrows=->,arrowsize=3pt 3.5](-2.9,2)(-1.6,2)
\psline[linecolor=blue,arrows=->,arrowsize=3pt 3.5](-2.9,1)(-2.1,1)

\rput(-1.5,3.3){\textrm{{\small $1$}}}
\rput(-1.2,2){\textrm{{\small $7$}}}
\rput(-2,0.7){\textrm{{\small $2$}}}
\rput(-0.7,1){\textrm{{\small $10$}}}
\rput(-1,-.3){\textrm{{\small $5$}}}

\psdots(4.5,2)(4.5,1)(4,0)(5,0)(5,-1)

{\psset{linecolor=blue} 
\psdots(2,2)(3,0)(3,1)(3,2)
\psdots(-6,2)(-5,2)(-4,2)(-3,2)(-3,1)}

\psline(4.5,3)(4.5,2)(4.5,1)(5,0)(5,-1)
\psline(4.5,1)(4,0)

\psline[linecolor=blue,arrows=->,arrowsize=3pt 3.5](2.1,2)(2.9,2)
\psline[linecolor=blue,arrows=->,arrowsize=3pt 3.5](3.1,2)(4.4,2)

\psline[linecolor=blue,arrows=->,arrowsize=3pt 3.5](3.1,0)(3.9,0)
\psline[linecolor=blue,arrows=->,arrowsize=3pt 3.5](3.1,1)(4.4,1)

\rput(4.5,3.3){\textrm{{\small $1$}}}
\rput(4.8,2){\textrm{{\small $7$}}}
\rput(4.8,1.1){\textrm{{\small $9$}}}
\rput(4.2,-.1){\textrm{{\small $2$}}}
\rput(5.3,-0.1){\textrm{{\small $10$}}}
\rput(5,-1.3){\textrm{{\small $5$}}}

\rput(2,2.3){\textrm{{\blue\small $3$}}}
\rput(3,2.3){\textrm{{\blue\small $4$}}}
\rput(2.8,1){\textrm{{\blue\small $6$}}}
\rput(2.8,0){\textrm{{\blue\small $8$}}}

\rput(0.8,1){\textrm{$\rightarrow$}}

\rput(-3,-1.2){\textrm{{\small lone-child vertices: 5,\, 7,}}}

\rput(-3,-1.7){\textrm{{\small and so $m=9$ and $c=7$}}}

\end{pspicture}
\end{center} 
\centerline{Figure 2b}

\vspace*{2mm}

After the conversion, $m$ has become a lone-child vertex in the tree and is still the max among the arc and lone-child vertices. The  conversion map is clearly reversible, and so defines an involution that changes the sign of the weight. Thus all weights cancel out except those  of the trees we wish to count, all of whose weights are 1, and the total weight, (\ref{eq1}), is the number of our trees.

\end{document}